\documentclass[12pt,a4paper]{article}
\usepackage[top=3.5cm, bottom=2.5cm, left=3.0cm, right=2.5cm]{geometry}
\usepackage{array,multirow,makecell}
\setcellgapes{1pt}
\makegapedcells
\newcolumntype{R}[1]{>{\raggedleft\arraybackslash }b{#1}}
\newcolumntype{L}[1]{>{\raggedright\arraybackslash }b{#1}}
\newcolumntype{C}[1]{>{\centering\arraybackslash }b{#1}}

\usepackage{color}
\usepackage{calc}
\usepackage[all]{xy}
\usepackage[centertags]{amsmath}
\usepackage{latexsym}
\usepackage{amsfonts}
\usepackage{amssymb}
\usepackage{amsthm}
\usepackage{amsfonts,graphicx}
\usepackage[colorlinks=true, linkcolor=blue, citecolor=blue, urlcolor=blue]{hyperref}
\numberwithin{equation}{section}
\usepackage[Glenn]{fncychap}
\usepackage{fancyhdr}
\usepackage{xcolor}
\usepackage{xcolor,rotating,epic,eepic}
\usepackage{tikz-qtree}
\usetikzlibrary{matrix}
\usepackage{fancyhdr}
\usepackage{xcolor,rotating,epic,eepic}
\usepackage{tikz}

\usepackage[babel=true,kerning=true]{microtype}

\usetikzlibrary{%
	arrows,%
	calc,%
	shapes.geometric,%
	shapes.misc,%
	shapes.symbols,%
	shapes.arrows,%
	automata,%
	through,%
	positioning,%
	scopes,%
	decorations.shapes,%
	decorations.text,%
	decorations.pathmorphing,%
	shadows}

\usepackage{fancyhdr}

\newenvironment{pv}{{\noindent {\bf Proof. }}}{\hfill {\rule{2mm}{2mm}} }

\newtheorem{thm}{Theorem}[section]
\newtheorem{cor}{Corollary}[section]
\newtheorem{lem}{Lemma}[section]
\newtheorem{prop}{Proposition}[section]

\newtheorem{dfn}{Definition}[section]

\theoremstyle{plaine}
\newtheorem{rem}{Remark}[section]

\usepackage{fancyhdr}
\title{ Dunkl symmetric coherent pairs of measures. An application to Fourier Dunkl-Sobolev expansions.}
\author{Mabrouk Sghaier\small$^1$\\ \footnotesize Higher Institute of Computer Medenine,   Medenine - 4119, Tunisia.\\
Francisco Marcell\'an, \small$^2$\\ \footnotesize Department of Mathematics, Universidad Carlos III de Madrid, 28911-Legan\'es, Spain\\
	Sabrine Hamdi\small$^3$  \\  \footnotesize Faculty of Sciences of Gabes, Laboratory of Mathematics and Applications (LR17ES11), 6072 Gabes, Tunisia.}
\date{}
\begin{document}
	\maketitle
	\textbf{Abstract}
	Let $\mathcal{T}_{\mu}$ be  the Dunkl operator. A pair of  symmetric measures $(u, v)$ supported on a symmetric subset of the real line is said to be a  symmetric Dunkl-coherent pair if the corresponding sequences of monic orthogonal polynomials $\{\mathrm{P}_n\}_{n\geq0}$ and $\{\mathrm{R}_n\}_{n\geq0}$ (resp.) satisfy
	\begin{equation*}
		\mathrm{R}_{n}(x)=\frac{\mathcal{T}_{\mu}\mathrm{P}_{n+1} (x)}{\mu_{n+1}}-\sigma_{n-1}\frac{\mathcal{T}_{\mu}\mathrm{P}_{n-1}(x)}{\mu_{n-1}}, \ \  n\geq2,
	\end{equation*}
	where $\{\sigma_n\}_{n\geq1}$ is a sequence of non-zero complex numbers and $\mu_{2n}=2n, \mu_{2n-1}= 2n-1+ 2\mu, n\geq1.$

 In this contribution we focus the attention on the sequence $\{\mathrm{S}_n^{(\lambda,\mu)}\}_{n\geq0}$ of monic orthogonal polynomials with respect to the Dunkl-Sobolev inner product
	\begin{equation*}
		<p,q>_{s,\mu}=<u,pq>+\lambda<v,\mathcal{T}_{\mu}p\mathcal{T}_{\mu}q>, \ \ \lambda >0, \ \ p, \ q  \ \in \mathcal{P}.
	\end{equation*}
	An algorithm is stated to compute the coefficients of the Fourier–Sobolev type  expansions with respect to $<. , .>_{s,\mu}$ for suitable smooth functions $f$ such that $f \in \mathcal{W}_2^1(\mathbb{R}, u, v, \mu)=\big\{ f; \   ||f||_{u}^{2} + \lambda || \mathcal{T}_{\mu }f||_{v} ^{2} <\infty\}$. Finally, two illustrative numerical  examples are presented.

\bigskip



	\section{Introduction}
	\ \  Let $\mathcal{P}$ be the linear space of polynomials with complex coefficients.  Let $u$ and $v$ denote two symmetric  measures supported on a symmetric subset of  the real line  and  let $\{\mathrm{P}_n\}_{n\geq0}$ and $\{\mathrm{R}_n\}_{n\geq0}$ be the corresponding sequences of monic orthogonal polynomials (MOPS in short). The pair $(u , v)$ is said to be a symmetric $\mathcal{T}_{\mu}$-coherent pair if there exists a sequence of non-zero complex numbers $\{\sigma_n\}_{n\geq1}$ such that
	\begin{equation*}
		\mathrm{R}_{n}(x)=\frac{\mathcal{T}_{\mu}\mathrm{P}_{n+1}(x)}{\mu_{n+1}}-\sigma_{n-1}\frac{\mathcal{T}_{\mu}\mathrm{P}_{n-1}(x)}{\mu_{n-1}}, \ \  n\geq2.
	\end{equation*}
	Here $\mathcal{T}_{\mu}$ is the Dunkl operator introduced  in  \cite{dunkl1991} that is defined as
	$$\mathcal{T}_{\mu}(f)= f'+2\mu\mathrm{H}_{-1}f, \ \ (\mathrm{H}_{-1}f)(x)=\frac{f(x)-f(-x)}{2x}, \ \ f \in \mathcal{P}, \ \mu\neq -n-\frac{1}{2},$$
and $\mu_{2n}=2n, \mu_{2n-1}= 2n-1+ 2\mu, n\geq1.$
	
Notice that in a more general framework the concept of symmetric $\mathcal{T}_{\mu}$-coherent pair of linear forms  is introduced   in \cite{sghaier2021}, where a description of them is stated.
	
In this paper, we consider the sequence  $\{\mathrm{S}_n^{(\lambda,\mu)}\}_{n\geq0}$ of  monic Dunkl-Sobolev polynomials orthogonal with respect to a special Sobolev inner product  ($\mathcal{T}_{\mu}$-Sobolev  inner product)
	\begin{equation*}
		<p,q>_{s,\mu}=<u,pq>+\lambda<v,\mathcal{T}_{\mu}p\mathcal{T}_{\mu}q>, \ \ \lambda >0, \ \ p, \ q  \ \in \mathcal{P}.
	\end{equation*}

Our results constitute an extension in the Dunkl framework of the paper \cite{Iserles1991}, where coherent and symmetrically coherent pairs of measures supported on the real line are introduced. Indeed,  a pair of positive and symmetric measures $(u, v)$ supported on a symmetric subset of  the real line is said to be a symmetrically coherent pair  if the corresponding MOPS  $\{P_{n}(x; u)\}_{n\geq0}, 	\{P_{n}(x; v)\}_{n\geq0} $ satisfy

$$	P_{n}(x; v) = \frac{1}{n+1}\left[P_{n+1}^{\prime}( x; u) -  \varrho_{n-1}  P_{n-1}^{\prime}( x; u)\right], \quad  \varrho_{n-1} \neq 0, \quad n \geq 2.$$

The pairs of symmetrically coherent positive measures  were completely determined  in \cite{Meijer1997}, where it is proved that one of the measures must be classical i. e., either Gegenbauer or  Hermite,  and the other one is a linear spectral transform (see \cite{ alexei}).  The motivations in \cite{Iserles1991} for introducing these coherent pairs of measures  are their applications in  Fourier analysis in terms of  sequences of orthogonal polynomials with respect to the  Sobolev inner products associated with the vector of measures $(u,v).$  In the Gegenbauer case, analytic properties of Sobolev-Gegenbauer polynomials were analyzed in \cite{Cleonice-10,  Oscar-SIGMA2018, maite, andrei1997}.  In the Hermite case, analytic properties of  Sobolev-Hermite  polynomials were analyzed in  \cite{Laura-03, Marcel-07, Paco-2006}. An exhaustive information can be found in the survey paper \cite{Paco2015}.\\

When we deal with symmetric $\mathcal{T}_{\mu}$-coherent pairs of measures, one of our main purposes is the effective computation of the Fourier coefficients in the expansions of  functions $f$ living in appropriate function spaces
$$f(x)\sim \sum_{n=0}^\infty \frac{f_n^{(\lambda, \mu)} }{s_n^{(\lambda,\mu)}}S_n^{(\lambda,\mu)}(x), $$ in terms of the  polynomials $\{\mathrm{S}_n^{(\lambda,\mu)}\}_{n\geq0}$  orthogonal with respect to the corresponding Sobolev inner product. This kind of Fourier expansions will be called  Dunkl Fourier-Sobolev expansions. Here  $s_n^{(\lambda,\mu)}=<\mathrm{S}_n^{(\lambda,\mu)},\mathrm{S}_n^{(\lambda,\mu)}>_{s,\mu}$ and $f_n^{(\lambda, \mu)} = <f,S_n^{(\lambda,\mu)}>_{s,\mu}$. \\

The structure of the paper is as follows. In Section $2$, we recall some basic definitions and results on the theory of orthogonal polynomials and the definition of the Dunkl operator, which will be used in the sequel. In Section $3$ we assume  $u$ (resp. $v$)  is a $\mathcal{T}_{\mu}$-classical symmetric measure, i. e., the family $\{\mathrm{P}_n\}_{n\geq0}$ (resp.  $\{\mathrm{R}_n\}_{n\geq0}$) is either the generalized  Hermite polynomial sequence or the generalized Gegenbauer polynomial sequence. The corresponding companion measure is also obtained in terms of either  a quadratic Chistoffel or a quadratic Geronimus transformation. In Section $4$, we consider Dunkl-Sobolev orthogonal polynomials  arising from  a symmetric $\mathcal{T}_{\mu}$-coherent pair of measures. In Section $5$ an  algorithm to compute Fourier Dunkl-Sobolev  coefficients is described. Finally, in Section $6$, a numerical example is discussed.

	\section{Preliminaries}
	\ \ Let $\mathcal{P}$ be the linear space of polynomials with complex coefficients and let $\mathcal{P}'$ be its algebraic dual space. We denote by $<u,f>$ the duality bracket for $u \in  \ \mathcal{P}'$ and $f \in \ \mathcal{P}$. In particular, we denote by $(u)_n:=<u,x^n>$, $n\geq0$, the moments of $u$. For a form $u$, any $a \in \mathbb{C}^*$,  $c \in \mathbb{C}$ and a polynomial $g$, let $Du=u'$, $gu$, $h_a u$, $\delta_c$ and $(x-c)^{-1}u$ be the forms defined, respectively,  by
	$$<u',f>:=-<u,f'>, \ <gu,f>:=<u,gf>, \
	<h_a u,f(x)>:=<u,h_af(x)> =<u,f(ax)> ,$$
	$$<\delta_c,f>:=f(c),\ <(x-c)^{-1}u,f(x)>:=<u,\theta_cf(x)>=<u,\frac{f(x)-f(c)}{x-c}>, \ \ f \in \mathcal{P}.$$
	 For $c,\,d \in \mathbb{C},\,c\neq d$, $u \in  \ \mathcal{P}'$ and $f \in \ \mathcal{P}$, it is s straightforward to prove that
	\begin{equation}\label{0}(x-c)^{-1}(x-c)u=u-(u)_0\delta_c.
\end{equation}

Let us recall some basic background on orthogonal polynomials and linear forms (see  \cite{chihara1978}).
	A form $u$ is called regular if there exists a sequence of polynomials $\{P_n\}_{n\geq0}, \deg P_n \leq n,$ such that
	$$\left\langle u , P_nP_m \right\rangle
	=p_n \delta_{nm},\quad p_n\neq0,\quad n\geq0\hspace*{0.1cm},$$
	where $\delta_{nm}$ is the Kronecker symbol. Then, $\deg
	P_n = n, n \geq 0,$ and we can always assume each $P_n$ is monic. In
	such a case, the sequence $\{P_n\}_{n\geq0}$ is unique. It is said
	to be the sequence of monic orthogonal polynomials (MOPS) with respect to
	$u$.
	
	There exist a complex sequence $\{\beta_n\}_{n\geq0}$ and a non zero complex sequence $\{\gamma_n\}_{n\geq1}$ such that the MOPS $\{P_n\}_{n\geq0}$ fulfils the following three-term recurrence relation (TTRR, in short)
	\begin{equation}\label{recurrence}\begin{array}{l}
			\ \ \ \ \ \ \ \ \ \ \ \  \ \ P_0(x)=1, \ \ \ P_1(x)=x-\beta_0,\\[0.2cm]
			P_{n+2}(x)=(x-\beta_{n+1})P_{n+1}(x)- \gamma_{n+1}P_{n}(x), \ n\geq0.
		\end{array}
	\end{equation}	
	By convention, we set $\gamma  _0 =(u)_0$. The form $u$ is said to be normalized if $(u)_{0}=1$. In this paper, we assume that any form will be normalized.
	
	A form $u$ is said to be  symmetric if and only if $(u)_{2n+1}=0$, $n\geq0$, or, equivalently, in \eqref{recurrence} $\beta_n=0$, $n\geq0$.
	\begin{prop} (\cite{chihara1978})
		Let $\{\mathrm{P}_n\}_{n\geq0}$ be the MOPS with respect to a
		regular form $u$. Then the following statements  are equivalent.
		\begin{itemize}
			
			\item[(i)] $u$ is symmetric.
			\item[(ii)] $\mathrm{P}_n(-x)= (-1)^n \mathrm{P}_n(x)$, $n\geq0$.
           \item[(iii)] There exist MOPS $\{\mathrm{A}_n\}_{n\geq0}, \{\mathrm{B}_n\}_{n\geq0},$ orthogonal with respect to linear forms $w, xw,$ respectively, such that $P_{2n}(x)= A_{n} (x^2), P_{2n+1}(x)= x B_{n} (x^2), n\geq0.$ On the other hand,  $w$ is the linear form with moments $\left\langle w ,  x^ {n}\right\rangle= \left\langle u ,  x^ {2n}\right\rangle.$
		\end{itemize}
	\end{prop}
	
	Let us introduce the Dunkl operator
	$$\mathcal{T}_{\mu}(f)= f'+2\mu\mathrm{H}_{-1}f, \ \ (\mathrm{H}_{-1}f)(x)=\frac{f(x)-f(-x)}{2x}, \ \ f \in \mathcal{P}, \ \mu\neq -n-\frac{1}{2}.$$
	This operator was introduced and studied for first time by Dunkl \cite{dunkl1991}. Notice that $\mathcal{T}_{0}$ is  the derivative operator $\mathrm{D}$. The transposed $t_{\mathcal{T}_{\mu}}$ of $\mathcal{T}_{\mu}$ is $t_{\mathcal{T}_{\mu}}=-\mathrm{D}-\mathrm{H}_{-1}=-\mathcal{T}_{\mu}.$ Thus, we have
	$$<\mathcal{T}_{\mu}u,f>=-<u,\mathcal{T}_{\mu}f>, \ \ u \in \mathcal{P}', \ \ f \in \mathcal{P}.$$
	In particular, this yields $<\mathcal{T}_{\mu}u,x^n>=-\mu_n<u,x^{n-1}>=-\mu_n(u)_{n-1}, \ n\geq0,$ where $(u)_{-1}=0$ and $\mu_n=n+\mu(1-(-1)^n),\ \ n \geq0$.\\[0.1cm]
	It is easy to see that
	\begin{equation}\label{du}
		\mathcal{T}_{\mu}(fu)=f\mathcal{T}_{\mu}u+f'u+2\mu(H_{-1}f)(h_{-1}u), \ \ f \in \mathcal{P}, \ \ u \in \mathbb{P'}.
	\end{equation}
	\begin{rem}
		Notice that if  $u$ is a symmetric form, then \eqref{du} becomes
		\begin{equation}\label{syu}
			\mathcal{T}_{\mu}(fu)=f\mathcal{T}_{\mu}u+\mathcal{T}_{\mu}fu, \ \ f \in \mathcal{P}.
		\end{equation}
	\end{rem}
	Now, let consider a sequence of monic polynomials $\{\mathrm{P}_n\}_{n\geq0}$ and let denote by
	\begin{equation}\label{eq1}
		\mathrm{P}^{[1]}_n(x,\mu)=\frac{\mathcal{T}_{\mu}\mathrm{P}_{n+1}(x)}{\mu_{n+1}}, \ \ \mu\neq -n-\frac{1}{2}
	\end{equation}
	a new sequence of monic polynomials.
	\begin{dfn}
		An MOPS $\{\mathrm{P}_n\}_{n\geq0}$ is said to be $\mathcal{T}_{\mu}$-classical if $\{\mathrm{P}^{[1]}_n(.,\mu)\}_{n\geq0}$ is also an MOPS. In this case, the associated form $u$ is called $\mathcal{T}_{\mu}$-classical form.
	\end{dfn}
	\begin{thm}( \cite{bencheikh2007})
		The only $\mathcal{T}_{\mu}$-classical symmetric MOPS are
		\begin{itemize}
	\item The generalized Hermite polynomials $\{\mathrm{H}_n^{\mu}\}_{n\geq0}$ for $\mu\neq -n-\frac{1}{2}$, $n\geq1$.
		
	\item The generalized Gegenbauer polynomials $\{\mathrm{G}_n^{(\alpha,\mu-\frac{1}{2})}\}_{n\geq0}$ for $\alpha\neq -n$, $\alpha+\mu \neq-n+\frac{1}{2}$,
		$\mu \neq -n + \frac{1}{2}$, $n\geq1$.
\end{itemize}
	\end{thm}
\begin{rem}

Generalized Hermite polynomials were introduced in \cite{ted}. Generalized Gegenbauer polynomials were studied in \cite{Belmehdi2001}. For more information about these two families the reader can check \cite{chihara1978}

\begin{prop}
	\begin{enumerate}
		\item The MOPS $\{\mathrm{G}_n^{(\alpha,\mu-\frac{1}{2})}\}_{n\geq0}$ satisfies the TTRR with
		\begin{gather*}
			\hat{\beta}_n=0,\ \ \ \gamma_{n+1}=\frac{(n+1+\varrho_n)(n+1+2\alpha+\varrho_n)}{4(n+\alpha+\mu+\frac{1}{2})(n+\alpha+\mu+\frac{3}{2})},\ \ n\geq0,\\
			\varrho_n=2\mu\frac{1+(-1)^n}{2}, \ \ n\geq0.
		\end{gather*}
	\item The MOPS $\{\mathrm{H}_n^{\mu}\}_{n\geq0}$ satisfies the TTRR with
	\begin{equation*}
		\hat{\beta}_n=0,\ \ \ \gamma_{n+1}=\frac{1}{2}(n+1+\mu(1+(-1)^n)),\ \ n\geq0.
	\end{equation*}	
	\end{enumerate}
\end{prop}

\end{rem}
	
In \cite{sghaier2012} a characterization of  $\mathcal{T}_{\mu}$-classical symmetric forms is given in terms of a distributional equation of Pearson type as
	\begin{thm}
		A symmetric form $u$ is said to be a $\mathcal{T}_{\mu}$-classical form if and only if $u$ is regular and there exist polynomials $\Phi$ and $\Psi , $  $\deg(\Phi)\leq2$ and $\deg(\Psi)=1,$ such that
		\begin{equation}\label{03}
			\mathcal{T}_{\mu}(\Phi u)=\Psi u,
		\end{equation}
		\begin{equation}\label{condition}
			\Psi'(0)-\frac{1}{2}\Phi''(0)\mu_n\neq0, \ \ n\geq0.
		\end{equation}
		The corresponding MOPS $\{\mathrm{P}_n\}_{n\geq0}$ is called  $\mathcal{T}_{\mu}$-classical.
	\end{thm}	

\begin{rem}

Notice that this is the $T_{\mu}$ analogue of the Pearson equation for the classical linear forms with respect to the derivative operator (see \cite{marcellan1994}).

\end{rem}

We state now the definition of the quasi-orthogonality.
\begin{dfn}
	Let $u \in \mathcal{P}'$ and $s$ be  a nonnegative integer number. A MPS $\{\mathrm{P}_n\}_{n\geq0}$ is said to be quasi-orthogonal of order $s$ with respect to $u$ if
	\begin{equation}\label{df1}
		\begin{array}{l}
			\big<u,P_mP_n\big>=0, \ \ 0\leq m \leq n-s-1, \ \ n\geq s+1,\\
			\ \ \ \ \ \ \ \ \ \ \ \ \ \ \ \ \exists \ r \geq s; \ \big<u,P_{r-s}P_r\big>\neq0.
		\end{array}
	\end{equation}
	
	If $\big<u,P_{r-s}P_r\big>\neq0$ for any $r\geq s$, then $\{\mathrm{P}_n\}_{n\geq0}$ is said to be strictly quasi-orthogonal of order $s$ with respect to $u$.
\end{dfn}

	\section{ Symmetric $\mathcal{T}_{\mu}$-Coherent Pairs}
	Symmetric $\mathcal{T}_{\mu}$-coherent pairs of linear forms have been introduced  by M. Sghaier and S. Hamdi in \cite{sghaier2021} as
	\begin{dfn}
		Let $u$ and $v$ denote two symmetric  forms and $\{\mathrm{P}_n\}_{n\geq0}$ and $\{\mathrm{R}_n\}_{n\geq0}$ will denote their respective SMOP. Assume that there exists a sequence of non-zero complex numbers $\{\sigma_n\}_{n\geq1}$ such that
		\begin{equation}\label{cohpair}
			\mathrm{R}_{n}(x)=\frac{\mathcal{T}_{\mu}\mathrm{P}_{n+1}(x)}{\mu_{n+1}}-\sigma_{n-1}\frac{\mathcal{T}_{\mu}\mathrm{P}_{n-1}(x)}{\mu_{n-1}}, \ \  n\geq2,
		\end{equation}
		holds. Then the pair $(u , v)$ is said to be a symmetric $\mathcal{T}_{\mu}$-coherent pair.
	\end{dfn}
	In \cite{sghaier2021} the symmetric $\mathcal{T}_{\mu}$-coherent pairs of linear forms are completely described and determined. Let us recall the important result.

	\begin{thm}
		If $(u, v)$ denote a symmetric $\mathcal{T}_{\mu}$-coherent pair of linear forms, then one of them is a symmetric Dunkl-classical form, i. e.,
		generalized Hermite or generalized Gegenbauer linear forms.
	\end{thm}

In the sequel we will assume  that $(u, v)$ is a symmetric $\mathcal{T}_{\mu}$-coherent pair of positive Borel measures. We will analyze two cases

\subsection{ $u$ is symmetric $\mathcal{T}_{\mu}$-classical }

 Let $u$ be a symmetric $\mathcal{T}_{\mu}$-classical measure, i. e.,  $\{\mathrm{P}_n\}_{n\geq0}$ is the sequence of generalized Hermite or generalized Gegenbauer polynomials. So, if
	$$\mathrm{T}_n(x)=\frac{\mathcal{T}_{\mu}\mathrm{P}_{n+1}(x)}{\mu_{n+1}}, \ \ n \geq0,$$
and taking into account that  $\{T_{n}\}_{ n\geq0}$ is again a $\mathcal{T}_{\mu}$-classical sequence of polynomials, of the same type, orthogonal with respect to a symmetric positive measure $w$, then  \eqref{cohpair} becomes
	\begin{equation}\label{coh2}
		\mathrm{R}_{n}(x)=\mathrm{T}_{n}(x)+\varepsilon_{n-2}\mathrm{T}_{n-2}(x), \ \  n\geq2,
	\end{equation}
	with $\varepsilon_{n-2}=-\sigma_{n-1}$ for all $n\geq2$. Furthermore, the coefficients of the TTRR  for the sequence $\{\mathrm{T}_n\}_{n\geq0}$, $\gamma_{n}$, are known.

\begin{prop}\label{prop}
 $\{R_{n}\}_{ n\geq0}$ is an MOPS if and only if
 	\begin{enumerate}
			\item [(i)] $\gamma_{n}+\varepsilon_{n-2} - \varepsilon_{n-1}>0$, $\forall \ n\geq1$.  By convention, $\varepsilon_{-1}=0.$
			\item [(ii)] $\varepsilon_n\gamma_{n}=(\gamma_{n+2} + \varepsilon_{n} - \varepsilon_{n+1}) \varepsilon_{n-1}, n \geq 1$.
		\end{enumerate}
\end{prop}

\begin{pv}
Since

\begin{equation}
				 \mathrm{R}_{n}(x)=x \mathrm{T}_{n}(x)+ \varepsilon_{n-2} \mathrm{T}_{n-2}(x),  ~~ n\geq1,
\end{equation}

from the TTRR that the MOPS $\{\mathrm{R}_n\}_{n\geq0}$ satisfies, you get

\begin{equation}				
x \mathrm{R}_{n}(x)=\mathrm{R}_{n+1}(x)+ (\gamma_{n} +\varepsilon_{n-2}  - \varepsilon_{n-1}) \mathrm{R}_{n-1}(x) +(\varepsilon_{n-2}\gamma_{n-2}- (\gamma_{n} + \varepsilon_{n-2} - \varepsilon_{n-1}) \varepsilon_{n-3})) T_{n-3}(x),  ~~n\geq1.
\end{equation}

Thus,  $\{\mathrm{R}_n\}_{n\geq0}$ satisfies a TTRR if and only if (i) and (ii) hold.
\end{pv}	

 In the next result, we compute the coefficients  $\{\tilde{\gamma}_{n} \}_{n\geq1}$ of
	the TTRR  for $\{\mathrm{R}_n\}_{n\geq0}$ in the terms of $\gamma_{n}$ and $\varepsilon_{n}$.

	\begin{cor}\label{relation}
		The following relations hold:
		\begin{enumerate}
			\item [(i)] $\tilde{\gamma}_{n}=\gamma_{n}+\varepsilon_{n-2} - \varepsilon_{n-1}$, $\forall \ n\geq1$.  By convention, $\varepsilon_{-1}=0.$
			\item [(ii)] $\varepsilon_n\gamma_{n}=\varepsilon_{n-1}\tilde{\gamma}_{n+2}$, $n \geq 1$.
		\end{enumerate}
	\end{cor}

	\begin{pv}
		Consider first the TTRR for the sequence $\{\mathrm{R}_n\}_{n\geq0}$:
		\begin{equation}\label{eq2}
			\begin{cases}
				\mathrm{R}_0(x)=1, \ \ 	\mathrm{R}_1(x)=x,\\
				\mathrm{R}_{n+1}(x)=x\mathrm{R}_{n}(x)-\tilde{\gamma}_{n} \mathrm{R}_{n-1}(x), \ \ n\geq1.
			\end{cases}
		\end{equation}
		Here $\tilde{\gamma}_{n} > 0$, $n \geq 1$.\\
		Substituting \eqref{coh2} into \eqref{eq2}, we get
		\begin{equation}\label{eq3}
			\begin{cases}
				\mathrm{T}_{n+1}(x)+\varepsilon_{n-1}\mathrm{T}_{n-1}(x)=x\big(\mathrm{T}_{n}(x)+\varepsilon_{n-2} \mathrm{T}_{n-2}(x)\big) -\tilde{\gamma}_{n}\big(\mathrm{T}_{n-1}(x)+\varepsilon_{n-3} \mathrm{T}_{n-3}(x)\big), \ \ \forall \ n\geq3,\\
				\mathrm{T}_{3}(x)+\varepsilon_1\mathrm{T}_{1}(x)=x\big(\mathrm{T}_{2}(x)+\varepsilon_{0}\big)  -\tilde{\gamma}_{2}\mathrm{T}_{1}(x),\\
				\mathrm{T}_{2}(x)+\varepsilon_0=x\mathrm{T}_{1} (x)-\tilde{\gamma}_{1},\\
				\mathrm{T}_{1}(x)=x.
			\end{cases}
		\end{equation}
		Consider now the TTRR  for the sequence $\{\mathrm{T}_n\}_{n\geq0}$:
		\begin{equation}\label{eq4}
			\begin{cases}
				\mathrm{T}_0(x)=1, \ \ 	\mathrm{T}_1(x)=x,\\
				x\mathrm{T}_{n}(x)=\mathrm{T}_{n+1}(x)+\gamma_{n} \mathrm{T}_{n-1}(x), \ \ n\geq1.
			\end{cases}
		\end{equation}
		Here $\gamma_{n} > 0$, $n \geq 1$.\\
		Using \eqref{eq4} and expanding $x\mathrm{T}_{i}(x)$ in \eqref{eq3} as a linear combination of the polynomials $\{\mathrm{T}_n\}_{n\geq0},$ we obtain
		\begin{equation}\label{eq5}
			\begin{cases}
				\mathrm{T}_{n+1}(x)+\varepsilon_{n-1}\mathrm{T}_{n-1}(x)=\mathrm{T}_{n+1}(x)+\gamma_{n} \mathrm{T}_{n-1}(x)+\varepsilon_{n-2}\mathrm{T}_{n-1}(x)+\varepsilon_{n-2}\gamma_{n-2} \mathrm{T}_{n-3}(x)\\
				\ \ \ \ \ \ \ \ \ \ \ \ \ \ \ \ \ \ \ \ \ \ \ \ \ \ \ \ \ \ \ \ \ \  -\tilde{\gamma}_{n}\big(\mathrm{T}_{n-1}(x)+\varepsilon_{n-3} \mathrm{T}_{n-3}(x)\big), \ \ \forall \ n\geq3,\\
				\mathrm{T}_{3}(x)+\varepsilon_1\mathrm{T}_{1}(x)=\mathrm{T}_{3}(x)+\gamma_{2} \mathrm{T}_{1}(x)+\varepsilon_{0}\mathrm{T}_{1} (x) -\tilde{\gamma}_{2}\mathrm{T}_{1}(x),\\
				\mathrm{T}_{2}(x)+\varepsilon_0=\mathrm{T}_{2}(x)+\gamma_{1} -\tilde{\gamma}_{1}.
			\end{cases}
		\end{equation}
		Since $\{\mathrm{T}_n\}_{n\geq0}$ is a basis in $\mathcal{P},$  after identification of coefficients in \eqref{eq5} we find
		\begin{equation}\label{1}
			\varepsilon_{n-1}=\gamma_{n}+\varepsilon_{n-2} -\tilde{\gamma}_{n}, \ \ \ \ n \geq3,
		\end{equation}
		\begin{equation}\label{2}
			\varepsilon_{1}=\gamma_{2}+\varepsilon_{0} -\tilde{\gamma}_{2},	
		\end{equation}
		\begin{equation}\label{3}
			\varepsilon_{0}=\gamma_{1}- \tilde{\gamma}_{1},
		\end{equation}
		\begin{equation}\label{4}
			\varepsilon_{n-2}\gamma_{n-2}-\varepsilon_{n-3}\tilde{\gamma}_{n}, \ \ n\geq3.
		\end{equation}
		
		From \eqref{1}-\eqref{3}, we get
		$$\varepsilon_{n-1}=\gamma_{n}+\varepsilon_{n-2}-\tilde{\gamma}_{n}, \ \ \ \ n \geq1, \  \ \varepsilon_{-1}=0. $$
		Thus  $(i)$ holds.  $(ii)$ can be deduced in a straightforward way from \eqref{4}.

	\end{pv}
	
	Proposition \ref{relation} requires the study of properties of the sequence $\{\varepsilon_{n} \}_{n\geq0}$. Next, we can prove that this sequence  is a solution of a non-linear difference equation.
	\begin{prop} The sequence $\{\varepsilon_{n} \}_{n\geq0}$ satisfies a quadratic difference equation
		\begin{equation}\label{eqdiff}
			\varepsilon_{n+1}+\frac{\gamma_{n}\gamma_{n+1}}{\varepsilon_{n-1}}=\gamma_{n+1}+\gamma_{n+2}+D,\ \ n \geq1,
		\end{equation}
		where
		$$D=(\varepsilon_{1}-\gamma_{2})(1-\frac{\gamma_{1}}{\varepsilon_{0}}).$$
	\end{prop}

	\begin{pv}
		From \ref{relation}, we have
		$$\varepsilon_n\gamma_{n}=\varepsilon_{n-1}(\gamma_{n+2}+\varepsilon_{n} - \varepsilon_{n+1}), \ \ n \geq 1.$$	
		As a consequence,
		$$\varepsilon_{n+1} - \varepsilon_{n}=\gamma_{n+2}-\frac{\varepsilon_n}{\varepsilon_{n-1}}\gamma_{n}, \ \ n \geq 1.$$
		We sum from $k=1$ to $k=n$ and we get
		\begin{equation}\label{sum1}
			\varepsilon_{n+1} - \varepsilon_{1} =\sum_{k=1}^{n} \gamma_{k+2} -\sum_{k=1}^{n} \frac{\varepsilon_k}{\varepsilon_{k-1}}\gamma_{k},  ~~n\geq 1.
		\end{equation}
		On the other hand, again from Corollary \ref{relation}, one can see that
		$$\frac{\varepsilon_{k-1}}{\varepsilon_{k-2}}\gamma_{k-1}=\gamma_{k+1}+\varepsilon_{k-1} - \varepsilon_{k}, \ \ k \geq 2.$$	
       Thus
		$$\frac{\varepsilon_{k}}{\varepsilon_{k-1}}=\frac{\gamma_{k+1}}{\varepsilon_{k-1}}-\frac{\gamma_{k-1}}{\varepsilon_{k-2}}+1, \ \ k \geq 2.$$
		Substitution of the previous equation into \eqref{sum1} yields
		\begin{align*}
			\varepsilon_{n+1} - \varepsilon_{1}&=\sum_{k=1}^{n}\gamma_{k+2}-\frac{\varepsilon_{1}\gamma_{1}}{\varepsilon_{0}}-\sum_{k=2}^{n}\Big(\frac{\gamma_{k}\gamma_{k+1}}{\varepsilon_{k-1}}-\frac{\gamma_{k-1}\gamma_{k}}{\varepsilon_{k-2}}+\gamma_{k}\Big)\\
			&= \gamma_{n+1}+\gamma_{n+2}-\gamma_{2}-\frac{\gamma_{n}\gamma_{n+1}}{\varepsilon_{n-1}} + \frac{\gamma_{1}\gamma_{2}}{\varepsilon_{0}}-\frac{\varepsilon_{1}\gamma_{1}}{\varepsilon_{0}}, \ \ n \geq 2.
		\end{align*}
		We denote by $D$ the term which does not depend on $n$ in the last relation, then
		$$D=\varepsilon_{1}-\frac{\varepsilon_{1}\gamma_{1}}{\varepsilon_{0}}+\frac{\gamma_{1}\gamma_{2}}{\varepsilon_{0}}-\gamma_{2}=(\varepsilon_{1}-\gamma_{2})(1-\frac{\gamma_{1}}{\varepsilon_{0}}),$$
		and we obtain \eqref{eqdiff}.
		
	\end{pv}

Notice that we can split the above difference equation for even and odd subindices, respectively. Indeed, for $n= 2m$ \eqref{eqdiff} reads
\begin{equation}
\varepsilon_{2m+1}+\frac{\gamma_{2m}\gamma_{2m+1}}{\varepsilon_{2m-1}}=\gamma_{2m+1}+\gamma_{2m+2}+D,\ \ m \geq1.
\end{equation}
 If we denote $\varepsilon_{2m+1}=\frac{ \rho_{m+1}}{\rho_{m}}, m\geq0, \rho_{0}=1,$  then the previous equation  becomes
 \begin{equation}
\rho_{m+1}+\gamma_{2m}\gamma_{2m+1} \rho_{m-1}=(\gamma_{2m+1}+\gamma_{2m+2}+D) \rho_{m},\ \ m \geq1.
\end{equation}
Here the initial conditions are  $\rho_{0}= 1, \rho_{1}= \varepsilon_{1}= \gamma D + \gamma-2,$ with $\gamma= (1- \frac{\gamma_{1}}{\varepsilon_{0}})^{-1}.$

Thus $\rho_{m}$ is a polynomial of degree $m$ in $D.$

For $n=2m-1,$ \eqref{eqdiff} reads
\begin{equation}
\varepsilon_{2m}+\frac{\gamma_{2m-1}\gamma_{2m}}{\varepsilon_{2m-2}}=\gamma_{2m}+\gamma_{2m+1}+D,\ \ m \geq1.
\end{equation}

 If we denote $\varepsilon_{2m}=\frac{ \tau_{m+1}}{\tau_{m}},  m\geq0, \tau_{0}=1,$ then the previous equation  becomes

 \begin{equation}
\tau_{m+1}+\gamma_{2m-1}\gamma_{2m} \tau_{m-1}= (\gamma_{2m}+\gamma_{2m+1}+D) \tau_{m},\ \ m \geq1.
\end{equation}

Here the initial conditions are  $\tau_{0}= 1, \tau_{1}= \varepsilon_{0}.$

In both cases, you have a TTRR.\\

	\begin{rem}
 Since the measure $u$ is known, we can compute the values of the coefficients $\varepsilon_{n}$ for $n\geq 2$, once the initial values  $\varepsilon_{0}$ and $\varepsilon_{1}$ are fixed.
 \end{rem}
 \bigskip
\begin{rem}
Notice that \eqref{coh2} means that the linear forms $w, v$ associated with the MOPS  $\{\mathrm{T}_n\}_{n\geq0}, \{\mathrm{R}_n\}_{n\geq0}$ are related by  $M_{2}(x) v= w,$ where $M_{2}$ is a positive quadratic even polynomial. In other words, $v$ is a Geronimus transformation of $w.$ When $u$ is the generalized Hermite form, then  $w=u$ and, as a consequence, $M_{2}(x) v= u,$  with $M_{2}(x)= \lambda(x^2 + \zeta ^2), \lambda >0.$ When $u$  is the generalized Gegenbauer  form, then $w= (1-x^{2})u$ and  $M_{2} (x) = \lambda (\zeta^2-x^2), |\zeta| \geq 1, \lambda >0,$ or $M_{2} (x) = \lambda (x^2+\zeta^2) , \zeta\in \mathbb{R}, \lambda >0.$
\end{rem}
\bigskip
An alternative way to get the sequence $\{\varepsilon_{n}\}_{n\geq 0}$ is the following one.\\

Let consider $T_{2n}(x)= A_{n}(x^2), T_{2n+1}(x)= x B_{n}(x^2), R_{2n}(x)= C_{n}(x^2), R_{2n+1}(x)=x D_{n}(x^2), n \geq 0.$ Thus \eqref{coh2}  reads
\begin{equation}
 C_{n}(x)= A_{n}(x) + \varepsilon_{2n-2} A_{n-1}(x), \ \ n\geq1,  D_{n}(x)= B_{n}(x) +\varepsilon_{2n-1} B_{n-1}(x),\ \ n\geq 1.
 \end{equation}

Taking into account that according to the Christoffel formula (see \cite{ted}) you get
\begin{equation}
xB_{n}(x)= A_{n+1}(x) + a_{n} A_{n}(x), xD_{n}(x)= C_{n+1}(x)+ b_{n} C_{n}(x), \ \ n \geq 0,
\end{equation}

then
\begin{equation}
A_{n+1}(x) + a_{n} A_{n}(x) + \varepsilon_{2n-1} ( A_{n}(x) + a_{n-1} A_{n-1}(x))= A_{n+1}(x) + \varepsilon_{2n} A_{n}(x)+ b_{n} (A_{n}(x) + \varepsilon_{2n-2} A_{n-1}(x)), \ \ n\geq1.
\end{equation}

As a consequence, identifying coefficients,
\begin{equation}
a_{n} + \varepsilon_{2n-1}= \varepsilon_{2n}+ b_{n}, \ \ n\geq 1,
\varepsilon_{2n-1} a_{n-1}= b_{n} \varepsilon_{2n-2}, \ \ n\geq 1.
\end{equation}

This is a coupled system for the sequence $\{ \varepsilon_{n}\}_{n\geq 0}$.\\

On the other hand, since $C_{n}(0)= A_{n}(0) + \varepsilon_{2n-2} A_{n-1}(0),\ \ n\geq1, ~~A_{n+1}(0) + a_{n} A_{n}(0)=0, \ \ n\geq0,$ we immediately get
\begin{equation}
a_{n}= -\frac{A_{n+1}(0)}{A_{n}(0)},  ~~ b_{n}= -\frac{C_{n+1}(0)}{C_{n}(0)}= -\frac{ A_{n+1}(0) + \varepsilon_{2n} A_{n}(0)} { A_{n}(0) + \varepsilon_{2n-2} A_{n-1}(0)}.
\end{equation}

Initial conditions for $\varepsilon_{0}, \varepsilon_{1}$ are needed  in order to generate the full sequence, Indeed,
\begin{equation}
\varepsilon_{2n}= a_{n} + \varepsilon_{2n-1} - \frac{\varepsilon_{2n-1} a_{n-1}}{\varepsilon_{2n-2}}.
\end{equation}
As a conclusion,
\begin{equation}
\varepsilon_{2n}=a_{n} + \varepsilon_{2n-1} (1- \frac{a_{n-1}}{\varepsilon_{2n-2}}), \ \ n\geq1.
\end{equation}

\begin{rem}
Notice that the MOPS $\{C_{n}\}_{n\geq0}$ is a canonical Geronimus transformation of the MOPS $\{A_{n}\}_{n\geq0}$ as well as  the MOPS $\{D_{n}\}_{n\geq0}$ is a canonical Geronimus transformation of the MOPS $\{B_{n}\}_{n\geq0},$ see \cite{pascal}.
\end{rem}

\subsection{ $v$ is symmetric $\mathcal{T}_{\mu}$-classical}

 Let $v$ be a symmetric $\mathcal{T}_{\mu}$-classical positive measure, i. e.,  $\{\mathrm{R}_n\}_{n\geq0}$ is the sequence of generalized Hermite or generalized Gegenbauer polynomials. So, if
	$$\mathrm{R}_n(x)=\frac{\mathcal{T}_{\mu}\mathrm{Q}_{n+1}(x)}{\mu_{n+1}}, \ \ n \geq0,$$
and taking into account that  $\{Q_{n}\}_{ n\geq0}$  is again a $\mathcal{T}_{\mu}$-classical sequence of polynomials,  orthogonal with respect to a symmetric  measure $\hat{v}$, then  \eqref{cohpair} becomes
	\begin{equation}\label{coh3}
		\mathrm{Q}_{n+1}(x)=\mathrm{P}_{n+1}(x)-\tau_{n}\mathrm{P}_{n-1}(x)- a_{n},\ \  n\geq1,
	\end{equation}
	with $\tau_{n}= \sigma_{n-1}\frac{\mu_{n+1}}{\mu_{n-1}},$ for all $n\geq2,$ and $\{a_{n}\}_{n\geq1}$ is a sequence of real numbers. Furthermore, the coefficients $\hat{\gamma}_{n}$ of the TTRR for the sequence $\{\mathrm{Q}_n\}_{n\geq0}$ are known.

Notice that $ a_{2n}= 0, n\geq0,$ as a straightforward consequence of the fact that $P_{2n-1}$ and $Q_{2n-1}$ are odd functions.

\begin{rem}
If $v$ is the generalized Hermite measure, then $\hat{v}= v.$ If $v$ is the generalized Gegenbauer measure, then $(1-x^2)\hat{v}= v.$
\end{rem}

\bigskip

Let consider $\mathrm{Q}_{2n}(x)= C_{n}(x^2), \mathrm{Q}_{2n+1}(x)= x D_{n}(x^2).$ Since  $\{\mathrm{Q}_{n}\}_{n\geq0}$ is an MOPS, then $\{C_{n}\}_{n\geq0}$ and   $\{D_{n}\}_{n\geq0}$   are MOPS and $x D_{n}(x)= C_{n+1}(x) + b_{n} C_{n}(x), \ \ n\geq0.$ On the other hand, $\mathrm{ P}_{2n}(x)= A_{n}(x^2), \mathrm{P}_{2n+1}=x B_{n}(x^2).$  Notice that  $\{\mathrm{P}_{n}\}_{n\geq0}$ is an MOPS  if and only if  $\{A_{n}\}_{n\geq0}$ and   $\{B_{n}\}_{n\geq0}$   are MOPS. Thus for $n=2m$ \eqref{coh3}  reads
\begin{equation}
 D_{m}(x)= B_{m}(x) - \tau_{2m} B_{m-1}(x), \ \ m\geq1.
 \end{equation}

Assuming $\tau_{2}\neq 0,$ if  $x D_{m}(x)= D_{m+1} (x) +\beta_{m} D_{m}(x) + \alpha_{m} D_{m-1}(x), \ \ m\geq1,$  and according to item 2, Theorem 1 in \cite{zecarlos}, then  $\{B_{m}\}_{m\geq0}$  is an MOPS if and only if $\tau_{2m}= \alpha_{m} \frac{ D_{m-1}(a)}{D_{m}(a)}, ~~ D_{m}(a)\neq 0, n\geq1, $ where $a= \beta_{0} + \frac{\alpha_{1}}{\tau_{2}}.$\\

If  $\{D_{n}\}_{n\geq0}$ is  the MOPS associated with the linear form  $\omega, $ then $\{B_{m}\}_{m\geq0}$  is an MOPS with respect to the linear form $\tilde{\omega}= (x-a) \omega,$ see \cite{zecarlos}.

\bigskip

On the other hand, for $n=2m-1$ \eqref{coh3}  reads
\begin{equation}
 C_{m}(x)= A_{m}(x) - \tau_{2m-1} A_{m-1}(x)- a_{2m-1}, \ \  m\geq1.
\end{equation}

Assuming $\tau_{1}+a_{1} \neq0,$ if  $x C_{m}(x)= C_{m+1} (x) +\tilde{\beta}_{m} C_{m}(x) + \tilde{\alpha}_{m} C_{m-1}(x), m\geq1,$  and according to item 2,  Theorem 1 in \cite{zecarlos},  then  $\{A_{m}\}_{m\geq0}$  is an MOPS if and only if $a_{2m-1}=0, \ \ m\geq1,$  and $\tau_{2m-1}= \tilde{\alpha}_{m} \frac{ C_{m-1}(b)}{C_{m}(b)}, ~~ C_{m}(b)\neq 0, \ \ n\geq1, $ where $b= \tilde{\beta}_{0} + \frac{\tilde{\alpha}_{1}}{\tau_{1}+ a_{1}}.$\\

If  $\{C_{n}\}_{n\geq0}$ is  the MOPS associated with the linear form  $\nu, $ then $\{A_{m}\}_{m\geq0}$  is an MOPS with respect to the linear form $\tilde{\nu}= (x-b) \nu,$ see item 2, Theorem 1 in \cite{zecarlos}.

\bigskip

Notice that $\tilde{\omega}= x\tilde{\nu},$ and  $\omega= x\nu.$ This means that $(x-a) x\nu= x\tilde{\nu}.$ But $\tilde{\nu}= (x-b) \nu.$ As a consequence, $ x (x-a) \nu= x (x-b) \nu.$ Thus $x(b-a) \nu=0. $ In other words, $(b-a) \nu= M \delta_{0}.$ Taking into account $\nu$ is regular, then $a=b.$
\bigskip

\begin{rem}

If $a_{1}= - \tau_{1},$ then  $\{A_{m}\}_{m\geq0}$  is an MOPS if and only if $a_{2m-1}=0$  and $\tau_{2m-1}=0, \ \ m \geq 2, $ a contradiction, see Theorem 1, item 1 in \cite{zecarlos}.

\end{rem}

\bigskip

Taking into account  that $\hat{v}_{2n}= \nu_{n}$ and $u_{2n}= \tilde{\nu}_{n},$ as a summary of the previous analysis, we get

\begin{prop}
If \eqref{coh3} holds, then $a_{n}=0 , n\geq 2.$ Moreover, $\langle \tilde{\nu}, x^{n} \rangle=\langle (x-b) \nu, x^{n} \rangle.$ This means that $u_{2n}=\hat{ v}_{2n+2}- b \hat{ v}_{2n}.$ In other words, $u=( x^2- b)\hat{ v},$ i. e., $u$ is a quadratic Christoffel transformation of the measure $\hat{ v}.$ In the generalized Hermite case, you have $(x^2 + \zeta^2) v= u,$ i. e., $b$ must be a nonpositive real number.  In the generalized Gegenbauer case,  $(1-x^2) u=( x^2+ b) v$ holds, i. e., either $b$ is a positive real number or   $(x^2-1) u=(x^2- b^2) v,$ where $|b| \geq 1.$
\end{prop}

The sequence $\{\tau_{n}\}_{n\geq 0}$ can be recursively obtained in the following way.

 According to the above notation  \eqref{coh2}  reads
\begin{equation}
 C_{n+1}(x)= A_{n+1}(x) - \tau_{2n+1} A_{n}(x),  \\ n\geq0,  D_{n}(x)= B_{n}(x) -\tau_{2n} B_{n-1}(x), \\ n\geq 1.
 \end{equation}

Taking into account that according to the Christoffel formula (see \cite{ted})  you get
\begin{equation}
xB_{n}(x)= A_{n+1}(x) + \hat{a}_{n} A_{n}(x), xD_{n}(x)= C_{n+1}(x)+ \hat{b}_{n} C_{n}(x), \\ n\geq 0,
\end{equation}

then
\begin{equation}
A_{n+1}(x) +\hat{ a}_{n} A_{n}(x) - \tau_{2n} ( A_{n}(x) + \hat{a}_{n-1} A_{n-1}(x))= A_{n+1}(x) - \tau_{2n+1} A_{n}(x)+ \hat{b}_{n} (A_{n}(x) - \tau_{2n-1} A_{n-1}(x)), \\ n\geq1.
\end{equation}

As a consequence, identifying coefficients,
\begin{equation}
\hat{a}_{n} -\tau_{2n}= -\tau_{2n+1}+ \hat{b}_{n}, \\ n\geq 1,
\tau_{2n} \hat{a}_{n-1}= \hat{b}_{n} \tau_{2n-1}, \\ n\geq 1.
\end{equation}

This is a coupled system for the sequence $\{ \tau_{n}\}_{n\geq 0}$. Initial conditions for $\tau_{1}$ are needed  in order to generate the full sequence,

\section{Dunkl-Sobolev Orthogonal Polynomials}

	Let $(u , v) $ be a $\mathcal{T}_{\mu}$-symmetric coherent pair of positive measures and let  $\{\mathrm{P}_n\}_{n\geq0}$ and $\{\mathrm{R}_n\}_{n\geq0}$  be the corresponding  MOPS , i. e.,  \eqref{cohpair} holds. We consider the Dunkl-Sobolev ($\mathcal{T}_{\mu}$-Sobolev) inner product as
	\begin{equation}\label{sobolev1}
		<p,q>_{s,\mu}=<u,pq>+\lambda<v,\mathcal{T}_{\mu}p\mathcal{T}_{\mu}q>, \ \ \lambda >0, \ \ p, \ q  \ \in \mathcal{P}.
	\end{equation}
	Let $\{\mathrm{S}_n^{(\lambda,\mu)}\}_{n\geq0}$ be the sequence of monic Dunkl-Sobolev polynomials orthogonal with respect to \eqref{sobolev1}. Notice that $\mathrm{P}_k=\mathrm{S}_k$ for $k=0,1,2$.
	
	Put
	$ <\mathrm{S}_n^{(\lambda,\mu)},\mathrm{S}_n^{(\lambda,\mu)}>_{s,\mu}=s_n^{(\lambda,\mu)}$, $p_n=<u,\mathrm{P}_n^2>$ and $r_n=<v,\mathrm{R}_n^2>$.
	
	The following theorem gives a fundamental algebraic property  for $\mathcal{T}_{\mu}$-symmetric coherent pairs of measures.

	\begin{thm}\label{thm1}

	Let $(u , v) $ be a  $\mathcal{T}_{\mu}$-symmetric coherent pair of positive measures  and let  $\{\mathrm{P}_n\}_{n\geq0}$ and $\{\mathrm{R}_n\}_{n\geq0}$ be the corresponding MOPS satisfying \eqref{cohpair}. Then  the Dunkl-Sobolev polynomials $\{S_n^{(\lambda,\mu)}\}_{n\geq0}$ are related to $\{\mathrm{P}_n\}_{n\geq0}$ and $\{\mathrm{R}_n\}_{n\geq0}$ by

		\begin{equation}\label{relation1}
			\mathrm{P}_{n+1}(x)+a_{n-2}\mathrm{P}_{n-1}(x)=S_{n+1}^{(\lambda,\mu)}(x)+\eta_{n-2}^{\mu} (\lambda)S_{n-1}^{(\lambda,\mu)}(x),\ \ n\geq2,
		\end{equation}
		\begin{equation}\label{relation2}
			\mu_{n+1}\mathrm{R}_n(x)=\mathcal{T}_{\mu}S_{n+1}^{(\lambda,\mu)}(x)+\eta_{n-2}^\mu(\lambda)\mathcal{T}_{\mu}S_{n-1}^{(\lambda,\mu)}(x),\ \ n\geq2,
		\end{equation}
		where
		\begin{equation}\label{coeff2}
			a_{n-2}=\frac{\mu_{n+1}}{\mu_{n-1}}\varepsilon_{n-2}, \ \ \eta_{n-2}^\mu(\lambda)=\frac{p_{n-1}}{s^{(\lambda,\mu)}_{n-1}}a_{n-2}, \ \ \ n \geq2.
		\end{equation}
		The coefficients $\eta_n^\mu(\lambda)$, $n\geq 0$, will be called Dunkl-Sobolev coefficients.
	\end{thm}
	\begin{pv}
		Since $\{\mathrm{S}_n^{(\lambda,\mu)}\}_{n\geq0}$ is a  basis in $\mathcal{P},$ then
		$$\frac{\mathrm{P}_{n+1}(x)}{\mu_{n+1}}-\sigma_{n-1}\frac{\mathrm{P}_{n-1}(x)}{\mu_{n-1}}=\frac{\mathrm{S}^{(\lambda,\mu)}_{n+1}(x)}{\mu_{n+1}}+\sum_{k=0}^{n}c_k^{(n)}\mathrm{S}^{(\lambda,\mu)}_k(x),$$
		for some constants $c_k^{(n)}.$ So, from \eqref{sobolev1} and \eqref{cohpair}, we get
		$$c_k^{(n)}\big<\mathrm{S}^{(\lambda,\mu)}_k,\mathrm{S}^{(\lambda,\mu)}_k\big>_{s,\mu}=\big<u,\Big(\frac{\mathrm{P}_{n+1}}{\mu_{n+1}}-\sigma_{n-1}\frac{\mathrm{P}_{n-1}}{\mu_{n-1}}\Big)\mathrm{S}^{(\lambda,\mu)}_k\big>+\lambda\big<v,\mathrm{R}_n\mathcal{T}_{\mu}\mathrm{S}^{(\lambda,\mu)}_k\big>.$$
		Thus,
		$$c_k^{(n)}=0, \ \ \forall \ 0 \leq k \leq n-2.$$
		So,
		$$\frac{\mathrm{P}_{n+1}(x)}{\mu_{n+1}}-\sigma_{n-1}\frac{\mathrm{P}_{n-1}(x)}{\mu_{n-1}}=\frac{\mathrm{S}^{(\lambda,\mu)}_{n+1}(x)}{\mu_{n+1}}+c_n^{(n)}\mathrm{S}^{(\lambda,\mu)}_{n} (x)+c_{n-1}^{(n)}\mathrm{S}^{(\lambda,\mu)}_{n-1}(x).$$
		Taking into account the fact that the polynomials $\mathrm{P}_{n+1}$ and $\mathrm{P}_{n-1}$ are both either even or odd functions, then $c_n^{(n)}=0.$
		Therefore,
		$$\frac{\mathrm{P}_{n+1}(x)}{\mu_{n+1}}-\sigma_{n-1}\frac{\mathrm{P}_{n-1}(x)}{\mu_{n-1}}=\frac{\mathrm{S}^{(\lambda,\mu)}_{n+1}(x)}{\mu_{n+1}}+c_{n-1}^{(n)}\mathrm{S}^{(\lambda,\mu)}_{n-1}(x).$$
		Furthermore,
		$$c_{n-1}^{(n)}=-\frac{\sigma_{n-1}}{\mu_{n-1}}\frac{p_{n-1}}{s_{n-1}^{(\lambda,\mu)}}.$$
		Thus \eqref{relation1} holds with
		$$a_{n-2}=\frac{\mu_{n+1}}{\mu_{n-1}}\varepsilon_{n-2}, \ \ \eta_{n-2}^\mu(\lambda)=\frac{p_{n-1}}{s_{n-1}^{(\lambda,\mu)}}a_{n-2}, \ \ n \geq 2.$$	
		Applying now the operator $\mathcal{T}_{\mu}$ to \eqref{relation1} and using \eqref{cohpair}, we obtain \eqref{relation2}.
	\end{pv}\\
	
	Notice that the converse result is not true. Indeed,
	\begin{thm}
		Let $u$ and $v$ be two symmetric positive measures and define the Dunkl-Sobolev inner product by \eqref{sobolev1}. We assume that the monic polynomials $\{\mathrm{P}_n\}_{n\geq0}$ are orthogonal with respect to  $u$ and let $\{\mathrm{S}_n^{(\lambda,\mu)}\}_{n\geq0}$ be  the sequence of Dunkl-Sobolev monic polynomials orthogonal with respect to \eqref{sobolev1}. Assume that they are related by \eqref{relation1}, with $\{a_n\}_{n\geq1}$ and $\{\eta_n^\mu(\lambda)\}_{n\geq1}$ nonzero sequences of real numbers. Then, there exists a sequence of real numbers $\{e_n\}_{n\geq2}$ such that
		$$\mathcal{T}_{\mu}\mathrm{P}_{n+1}(x)+a_{n-2}\mathcal{T}_{\mu}\mathrm{P}_{n-1}(x)=\mathrm{R}_{n}(x)+e_n \mathrm{R}_{n-2}(x), \ \ n \geq 2,$$
		where $\{\mathrm{R}_n\}_{n\geq0}$ denotes the MOPS with respect to the measure $v$.		
	\end{thm}
	\begin{pv}
		Notice that
		$$\mathcal{T}_{\mu}\mathrm{P}_{n+1}(x)+a_{n-2}\mathcal{T}_{\mu }\mathrm{P}_{n-1}(x)=\mu_{n+1}\mathrm{R}_n(x)+\sum_{i=0}^{n-1}e_{n,k}\mathrm{R}_k(x),$$
		where
		$$e_{n,k}=\frac{<v,\big(\mathcal{T}_{\mu}\mathrm{P}_{n+1}+a_{n-2}\mathcal{T}_{\mu }\mathrm{P}_{n-1}\big)\mathrm{R}_k>}{r_k}.$$
		On the other hand, from \eqref{sobolev1} and \eqref{relation1}, for $0\leq k \leq n-2$, we get
		$$0=<S_{n+1}^{(\lambda,\mu)}+\eta_{n-2}^\mu(\lambda)S_{n-1}^{(\lambda,\mu)},x^k>_{s,\mu}=\lambda \mu_k<v,\big(\mathcal{T}_{\mu}\mathrm{P}_{n+1}+a_{n-2}\mathcal{T}_{\mu }\mathrm{P}_{n-1}\big) x^{k-1}>.$$
		Then,
		$$<v,\big(\mathcal{T}_{\mu}\mathrm{P}_{n+1}+a_{n-2}\mathcal{T}_{\mu }\mathrm{P}_{n-1}\big) x^{k-1}>=0, \ \  for \ all \ \ 0\leq k \leq n-3.$$
		So, we can write $\mathcal{T}_{\mu}\mathrm{P}_{n+1}+a_n\mathcal{T}_{\mu }\mathrm{P}_{n-1}$ as
		$$\mathcal{T}_{\mu}\mathrm{P}_{n+1}+a_{n-2}\mathcal{T}_{\mu }\mathrm{P}_{n-1}=\mu_{n+1}\mathrm{R}_n+e_{n,n-1}\mathrm{R}_{n-1}+e_{n,n-2}\mathrm{R}_{n-2}, \ \ n\geq2.$$
		Because of the symmetry of the polynomials, the coefficient $e_{n,n-1}$ most be zero. Let denote $e_{n,n-2}=e_n$. Then
		$$\mathcal{T}_{\mu}\mathrm{P}_{n+1 (x)}+a_{n-2}\mathcal{T}_{\mu }\mathrm{P}_{n-1}(x)=\mu_{n+1}\mathrm{R}_{n}(x)+e_n\mathrm{R}_{n-2}(x), \ \ n\geq2.$$
		Now we give an explicit expression for the coefficients $e_n$. Using the orthogonality of the sequence  $\{\mathrm{S}_n^{(\lambda,\mu)}\}_{n\geq0}$ with respect to \eqref{sobolev1}, we have
		$$<S_{n+1}^{(\lambda,\mu)}+\eta_{n-2}^\mu(\lambda)S_{n-1}^{(\lambda,\mu)},x^{n-1}>_{s,\mu}=\eta_{n-2}^\mu(\lambda)s_{n-1}^{(\lambda,\mu)}.$$
		On the other hand, by \eqref{relation1}, we get
		$$<S_{n+1}^{(\lambda,\mu)}+\eta_{n-2}^\mu(\lambda)S_{n-1}^{(\lambda,\mu)},x^{n-1}>_{s,\mu}=a_{n-2}p_{n-1}+\lambda\mu_{n-1}e_nr_{n-2}, \ \  n\geq 2.$$
		Therefore,
		$$e_n= \frac{\eta_{n-2}^\mu(\lambda)s_{n-1}^{(\lambda,\mu)}-a_{n-2}p_{n-1}}{\lambda \mu_{n-1}r_{n-2}},\ \ n\geq2.$$
	\end{pv}\\

\begin{rem}
This is the analogue in the Dunkl framework of a result stated in \cite{Delgado2005}.
\end{rem}
	
	From the previous theorem, we see that $(u , v)$ is a symmetric $\mathcal{T}_{\mu}$-coherent pair of positive measures  if and only if
	$e_n=0$ for all $n\geq2$, which is equivalent to
	$$ \eta_{n-2}^\mu(\lambda)=\frac{p_{n-1}}{s_{n-1}^{(\lambda,\mu)}}a_{n-2}, \ \ n \geq 2.$$
	
	The following lemma provides a  recurrence relation for the sequence $\{\mathrm{s}_n^{(\lambda,\mu)}\}_{n\geq0}$.
	\begin{lem}
		For $n \geq 2$
		\begin{equation}\label{sobl1}
			s_{n+1}^{(\lambda,\mu)}=p_{n+1} + a_{n-2}^2p_{n-1}+   \lambda \mu_{n+1}^2 r_{n}-(\eta_{n-2}^\mu(\lambda))^2 s_{n-1}^{(\lambda,\mu)}.
		\end{equation}
	\end{lem}
	\begin{pv}
		For $n\geq2$, we define
		$$\mathrm{Q}_{n+1}(x)=\mathrm{P}_{n+1}(x)+a_{n-2}\mathrm{P}_{n-1}(x).$$
		Then, from \eqref{relation1},
		\begin{equation}\label{Q}
			\mathrm{Q}_{n+1}(x)=S_{n+1}^{(\lambda,\mu)}(x)+\eta_{n-2}^\mu(\lambda)S_{n-1}^{(\lambda,\mu)}(x).
		\end{equation}
		We first prove  that for all $n \geq 4$
		\begin{equation}\label{Qnu}
			\eta_{n-2}^\mu(\lambda)=\frac{<\mathrm{Q}_{n+1},\mathrm{Q}_{n-1}>_{s,\mu}}{<\mathrm{Q}_{n-1}, \mathrm{Q}_{n-1}>_{s,\mu}-\eta_{n-4}^\mu(\lambda)<\mathrm{Q}_{n-1},\mathrm{Q}_{n-3}>_{s,\mu}}.
		\end{equation}
As  a straightforward consequence of \eqref{Q} we get
		\begin{equation}\label{Q1}
			<\mathrm{Q}_{n+1},\mathrm{Q}_{n-1}>_{s,\mu}=\eta_{n-2}^\mu(\lambda)<S_{n-1}^{(\lambda,\mu)},\mathrm{Q}_{n-1}>_{s,\mu}
		\end{equation}
		and
		\begin{align}\label{Q2}
			\nonumber
			<\mathrm{Q}&_{n-1},\mathrm{Q}_{n-1}-\eta_{n-4}^\mu(\lambda)\mathrm{Q}_{n-3}>_{s,\mu}\\\nonumber
			&=<S_{n-1}^{(\lambda,\mu)}+\eta_{n-4}^\mu(\lambda)S_{n-3}^{(\lambda,\mu)},\mathrm{Q}_{n-1}-\eta_{n-4}^\mu(\lambda)\mathrm{Q}_{n-3}>_{s,\mu}\\\nonumber
			&=<S_{n-1}^{(\lambda,\mu)},\mathrm{Q}_{n-1}>_{s,\mu}+\eta_{n-4}^\mu(\lambda) <S_{n-3}^{(\lambda,\mu)},\mathrm{Q}_{n-1}>_{s,\mu}-\big(\eta_{n-4}^\mu(\lambda)\big)^2<S_{n-3}^{(\lambda,\mu)},\mathrm{Q}_{n-3}>_{s,\mu}\\\nonumber
			&=<S_{n-1}^{(\lambda,\mu)},\mathrm{Q}_{n-1}>_{s,\mu}+\eta_{n-4}^\mu(\lambda) <S_{n-3}^{(\lambda,\mu)},S_{n-1}^{(\lambda,\mu)}+\eta_{n-4}^\mu(\lambda)S_{n-3}^{(\lambda,\mu)}>_{s,\mu}-\big(\eta_{n-4}^\mu(\lambda)\big)^2<S_{n-3}^{(\lambda,\mu)},\mathrm{Q}_{n-3}>_{s,\mu}\\\nonumber
			&=<S_{n-1}^{(\lambda,\mu)},\mathrm{Q}_{n-1}>_{s,\mu} +\big(\eta_{n-4}^\mu(\lambda)\big)^2<S_{n-3}^{(\lambda,\mu)},\mathrm{Q}_{n-3}>_{s,\mu}-\big(\eta_{n-4}^\mu(\lambda)\big)^2<S_{n-3}^{(\lambda,\mu)},\mathrm{Q}_{n-3}>_{s,\mu}\\
			&=<S_{n-1}^{(\lambda,\mu)},\mathrm{Q}_{n-1}>_{s,\mu}.
		\end{align}
		Replacing \eqref{Q2} in \eqref{Q1}
		$$<\mathrm{Q}_{n+1},\mathrm{Q}_{n-1}>_{s,\mu}=\eta_{n-2}^\mu(\lambda)<\mathrm{Q}_{n-1},\mathrm{Q}_{n-1}-\eta_{n-4}^\mu(\lambda)\mathrm{Q}_{n-3}>_{s,\mu}$$
		and \eqref{Qnu} holds.
		
		Now, from the definition,
		\begin{align*}
			<\mathrm{Q}_{n-1},\mathrm{Q}_{n-1}>_{s,\mu}&=<u,\mathrm{Q}_{n-1}\mathrm{Q}_{n-1}>+\lambda<v,\mathcal{T}_{\mu}\mathrm{Q}_{n-1}\mathcal{T}_{\mu}\mathrm{Q}_{n-1}>\\
			&=<u,\mathrm{Q}_{n-1}\mathrm{Q}_{n-1}>+\lambda \mu_{n-1}^2<v,\mathrm{R}_{n-2}(x)\mathrm{R}_{n-2}(x)>\\
			&=<u,(\mathrm{P}_{n-1}+a_{n-4}\mathrm{P}_{n-3})(\mathrm{P}_{n-1}+a_{n-4}\mathrm{P}_{n-3})>+\lambda \mu_{n-1}^2 r_{n-2}\\
			&=p_{n-1}+a^2_{n-4}p_{n-3}+\lambda \mu_{n-1}^2 r_{n-2},
		\end{align*}
		\begin{align*}
			<\mathrm{Q}_{n+1},\mathrm{Q}_{n-1}>_{s,\mu}&=<u,\mathrm{Q}_{n+1}\mathrm{Q}_{n-1}>+\lambda<v,\mathcal{T}_{\mu}\mathrm{Q}_{n+1}\mathcal{T}_{\mu}\mathrm{Q}_{n-1}>\\
			&=<u,\mathrm{Q}_{n+1}\mathrm{Q}_{n-1}>+\lambda\mu_{n+1}\mu_{n-1}<v,\mathrm{R}_{n}(x)\mathrm{R}_{n-2}(x)>\\
			&=<u,(\mathrm{P}_{n+1}+a_{n-2}\mathrm{P}_{n-1})(\mathrm{P}_{n-1}+a_{n-4}\mathrm{P}_{n-3})>\\
			&=a_{n-2} p_{n-1}
		\end{align*}
		and
		$$<\mathrm{Q}_{n-1},\mathrm{Q}_{n-3}>_{s,\mu}=a_{n-4} p_{n-3}.$$
		So, \eqref{Qnu} becomes
		\begin{equation}\label{Qnu2}
			\eta_{n-2}^\mu(\lambda)=\frac{a_{n-2} p_{n-1}}{p_{n-1}+a^2_{n-4}p_{n-3}+\lambda \mu_{n-1}^2 r_{n-2}-\eta_{n-4}^\mu(\lambda)a_{n-4} p_{n-3}}.
		\end{equation}
		According to \eqref{coeff2}, \eqref{sobl1} holds.
	\end{pv}
	\begin{rem}
	\eqref{sobl1} is useful in order to compute the norms  $s_{n}^{(\lambda,\mu)}$ if the Dunkl-Sobolev coefficients are known, with some initial conditions.
		First, from their definitions it is easy to see that
		\begin{equation}\label{algo1}
			s_{1}^{(\lambda,\mu)}=\lambda + p_1, \ \ \ \eta_{0}^\mu(\lambda)=\frac{a_0 p_1}{\lambda + p_1}.
		\end{equation}
		Now, for $n=1$,
		$$\eta_{1}^\mu(\lambda)=\frac{a_1p_2}{s_{2}^{(\lambda,\mu)}}.$$
		Since $S_2^{(\lambda,\mu)}(x)=\mathrm{P}_{2}(x)= x^2-\gamma_{1},$ then
		\begin{equation}\label{algo2}
			s_{2}^{(\lambda,\mu)}=p_2+4\lambda r_1, \ \ \ \eta_{1}^\mu(\lambda)=\frac{a_1p_2}{p_2+4\lambda r_1}.
		\end{equation}
		
		Now, using \eqref{algo1} and from \eqref{sobl1} with $n=2$ we obtain $s_{3}^{(\lambda,\mu)}$. As a consequence, from \eqref{coeff2} we get $\eta_{2}^\mu(\lambda)$.
		
		In an analogue way, for $n=2k$, where $k$ is an integer number greater than or equal to $2$,  in \eqref{sobl1} we obtain $s_{2k+1}^{(\lambda,\mu)}$ and then $\eta_{2k}^\mu(\lambda)$.
		
		Similarly, using \eqref{algo2} and from \eqref{sobl1} with $n=2k+1$, where $k$ is an integer number greater than or equal to $2$, we can get $s_{2k+2}^{(\lambda,\mu)}$ and then $\eta_{2k+1}^\mu(\lambda)$.
	\end{rem}
	
	\section{Algorithm for  the Fourier Dunkl-Sobolev coefficients}

	Let us define  the Dunkl-Sobolev space as

	$$\mathcal{W}_2^1(\mathbb{R}, u, v, \mu)=\big\{ f;  \   ||f||^{2}_{u} +\lambda ||\mathcal{T}_{\mu} f ||^{2}_{v} <\infty\}.$$
	
	If  $f \in \mathcal{W}_2^1(\mathbb{R}, u, v, \mu)$, then its Fourier expansion in terms of Dunkl-Sobolev  orthogonal polynomials $\{S_n^{\lambda_\mu}\}_{n\geq0}$ is

	$$f(x) \sim \sum_{n=0}^\infty \frac{<f,S_n^{(\lambda,\mu)}>_{s,\mu}}{s_n^{(\lambda,\mu)}}S_n^{(\lambda,\mu)}(x).$$

	We denote $f_n^{(\lambda, \mu)}= <f,S_n^{(\lambda,\mu)}>_{s,\mu}$ and $\mathcal{F}^{(\lambda, \mu)}_n=\frac{f_n^{(\lambda, \mu)}}{s_n^{(\lambda,\mu)}}$. $\mathcal{F}^{(\lambda, \mu)}_n$ is said to be the $n$-th Dunkl Fourier-Sobolev coefficient.
	
	Now, take $f \in \mathcal{W}_2^1(\mathbb{R}, u, v, \mu)$ and let consider a sequence of real numbers $\{w_n^\mu(f)\}_{n\geq0}$ such that
	\begin{equation}\label{w n mu}
		w_n^\mu(f)=<u,(\mathrm{P}_{n+2}+a_{n-1}\mathrm{P}_n)f>+\lambda<v,(\mathcal{T}_{\mu}\mathrm{P}_{n+2}+ a_{n-1}\mathcal{T}_{\mu}\mathrm{P}_n)\mathcal{T}_{\mu}f>.
	\end{equation}
	
	\begin{lem}
		The following relation holds
		\begin{equation}\label{fourier}
			f_{n+2}^{(\lambda, \mu)}+ \eta_{n-1}^\mu(\lambda) f_n^{(\lambda, \mu)}= w_n^\mu(f), \ \ \ n\geq0,
		\end{equation}
		with  initial conditions
		$$\eta_{-1}^\mu(\lambda)=0,\ \  f_0^{(\lambda, \mu)}=<f,1>_{s,\mu}=<u,f>, \ \ f_1^{(\lambda, \mu)}=<u,xf>+\lambda \mu_1 <v,\mathcal{T}_{\mu}f>$$
		and
		$$w_0^\mu(f)=<u,\mathrm{P}_{2}f>+\lambda<v,\mathcal{T}_{\mu}\mathrm{P}_{2}\mathcal{T}_{\mu}f>.$$
	\end{lem}
	\begin{pv}
		Using \eqref{relation1}, we get
		\begin{align*}
			<f,S_{n+2}^{(\lambda,\mu)}>_{s,\mu} &= -\eta_{n-1}^\mu(\lambda)<f,S_n^{(\lambda,\mu)}>_{s,\mu}+<f,(\mathrm{P}_{n+2}+ a_{n-1}\mathrm{P}_n)>_{s,\mu}\\
			&=-\eta_{n-1}^\mu(\lambda)<f,S_n^{(\lambda,\mu)}>_{s,\mu}+<u,(\mathrm{P}_{n+2}+a_{n-1}\mathrm{P}_n)f>\\ & \ \ \ \ \ \ \ \ \ \ \ \ \ \ \ \ \ \ \ \ \ \ \ \ \ \ \  \ \ \ \ \ \ +\lambda<v,(\mathcal{T}_{\mu}\mathrm{P}_{n+2}+
			a_{n-1}\mathcal{T}_{\mu}\mathrm{P}_n)\mathcal{T}_{\mu}f>
		\end{align*}
		and the result follows.
	\end{pv}
	
	We will now describe an algorithm to compute the Fourier  Dunkl-Sobolev coefficients  in the Fourier expansions in terms of  Dunkl-Sobolev orthogonal polynomials. We assume that $\varepsilon_{0}$ and $\varepsilon_{1}$ are known.\\
	
	\textbf{Algorithm 1.}\textbf{\emph{ Even order  Fourier Dunkl-Sobolev coefficients}}\\
	
	 For $n$ even, the  Fourier Dunkl-Sobolev coefficients $\mathcal{F}^{(\lambda, \mu)}_n$ can be computed as follows.\\
	
	 Step $0$. Initial conditions $\lambda$, $\mu$, $\varepsilon_{0}$, $\varepsilon_{1}$, $\eta_{-1}^\mu(\lambda)$, $s_0^{(\lambda,\mu)}$, $f_0^{(\lambda, \mu)}$,  $w_0^\mu(f)$, $s_2^{(\lambda,\mu)}$ and  $\eta_{1}^\mu(\lambda)$.\\
	
	Step $1$. Using step $0$ and from \eqref{fourier} when $n=0$ you get  $f_2^{{(\lambda, \mu)}}.$  Then calculate $F_2^{(\lambda, \mu)}$.\\
	
	Step $2$. Using step $0$ and the information in step $1$ to compute: $w_2^\mu(f)$ with \eqref{w n mu} and $n=2$, $f_4^{(\lambda, \mu)}$ from \eqref{fourier} with $n=2$ and $s_4^{(\lambda,\mu)}$ taking $n=3$ in \eqref{sobl1}. Then compute $F_4^{(\lambda, \mu)}$.\\
	
	Step $k$. For $k\geq 3$, using the step $0$ and the information in steps $1$ to $k-1$ to compute: $\varepsilon_{2k-3}$ from \eqref{eqdiff} and $n=2k-4$,  $a_{2k-3}$ from \eqref{coeff2} with $n=2k-1$, $\eta_{2k-3}^\mu(\lambda)$ from \eqref{coeff2} with $n=2k-1$, $w_{2k-2}^\mu(f)$ with \eqref{w n mu} and $n=2k-2$, $f_{2k}^{(\lambda, \mu)}$ from \eqref{fourier} with $n=2k-2$ and, finally, $s_{2k}^{(\lambda,\mu)}$ taking $n=2k-1$ in \eqref{sobl1}. Then, compute $F_{2k}^{(\lambda, \mu)}$.\\

	\textbf{Algorithm 2.} \textbf{\emph{Odd order  Fourier Dunkl-Sobolev coefficients}}\\
	
	For $n$ odd, the Fourier Dunkl-Sobolev coefficients $\mathcal{F}^{(\lambda, \mu)}_n$ can be computed as follows.\\
	
	Step $0$. Initial conditions $\lambda$, $\mu$, $\varepsilon_{0}$, $\varepsilon_{1}$, $\eta_{0}^\mu(\lambda)$, $s_1^{(\lambda,\mu)}$, $f_1^{(\lambda, \mu)}$ and $w_1^\mu(f)$.\\
	
	Step $1$. Using step $0$  to compute  $f_2^{{(\lambda, \mu)}}$ from \eqref{fourier} with $n=0$  and $s_3^{(\lambda,\mu)}$ from \eqref{sobl1} with $n=2$. Then, calculate $F_3^{(\lambda, \mu)}$.\\
	
	Step $2$. Using step $0$ and the information in step $1$ to compute: $\varepsilon_{2}$ from \eqref{eqdiff} and $n=1$, $a_{2}$ from \eqref{coeff2} with $n=4$, $\eta_{2}^\mu(\lambda)$ from \eqref{coeff2} with $n=4$, $w_3^\mu(f)$ with \eqref{w n mu} and $n=3$, $f_5^{(\lambda, \mu)}$ from \eqref{fourier} with $n=3$ and $s_5^{(\lambda,\mu)}$ taking $n=4$ in \eqref{sobl1}. Then compute $F_5^{(\lambda, \mu)}$.\\
	
	Step $k$. For $k\geq 3$, using the step $0$ and the information in steps $1$ to $k-1$ to compute: $\varepsilon_{2k-2}$ from \eqref{eqdiff} and $n=2k-3$,  $a_{2k-2}$ from \eqref{coeff2} with $n=2k$, $\eta_{2k-2}^\mu(\lambda)$ from \eqref{coeff2} with $n=2k$, $w_{2k-1}^\mu(f)$ with \eqref{w n mu} and $n=2k-1$, $f_{2k+1}^{(\lambda, \mu)}$ from \eqref{fourier} with $n=2k-1$ and finally $s_{2k+1}^{(\lambda,\mu)}$ taking $n=2k$ in \eqref{sobl1}. Then, compute $F_{2k+1}^{(\lambda, \mu)}$.\\

\begin{rem}
Notice that an algorithm for generate Fourier-Sobolev coefficients when you deal with $(1,1)$ coherent pairs of measures with respect to the derivative operator  has been presented in \cite{herbert}.
\end{rem}

	\section{Numerical examples}
	In order to illustrate the algorithm presented above, in this section we give an example involving the construction of the  Fourier Dunkl–Sobolev series with respect to $\mathcal{T}_{\mu}-$ Sobolev inner products of the type \eqref{sobolev1} where $(u, v)$ is  a $\mathcal{T}_{\mu}$-symmetric coherent pair of measures and $u$ is the classical one, i.e. $u$ is the generalized Hermite weight  $\mathcal{H}^{\mu}$  or the generalized Gegenbauer weight $\mathcal{G}^{(\alpha,\mu-\frac{1}{2})}$.
	\subsection{The generalized Hermite case}
	The sequence of generalized Hermite polynomials $\{\mathrm{H}_n^{\mu}\}_{n\geq0}$ is an MOPS with respect to the positive definite linear form $\mathcal{H^\mu}$ defined by the weight function
	$$\omega(x)=\frac{1}{\Gamma(\alpha+\frac{1}{2})}|x|^{2\mu}e^{-x^2}, \mu>-1/2,$$
supported on the real line.

	This positive measure  is $\mathcal{T}_{\mu}$-classical and satisfies the $\mathcal{T}_{\mu}$-Pearson equation
 $$\mathcal{T}_{\mu}(\mathcal{H^\mu})=-2x\mathcal{H^\mu}.$$

	If we apply $\mathcal{T}_{\mu}$ to $\mathcal{H}_n^{\mu}$, we obtain \cite{Ros1994}
	\begin{equation}\label{hermet}
		\mathcal{T}_{\mu}\mathrm{H}_n^{\mu}(x)=\mu_n\mathrm{H}_{n-1}^{\mu}(x).
	\end{equation}

This means that $\{\mathrm{H}_n^{\mu}\}_{n\geq0}$ is an $\mathcal{T}_{\mu}$-Appell sequence.

Now, let $\{u , v\}$ be a $\mathcal{T}_{\mu}$-symmetrically coherent pair of positive measures. Suppose that  $u$ is $\mathcal{T}_{\mu}$-classical with $u=\mathcal{H}^{\mu}$. From \cite{sghaier2021}, we obtain
	$$\mathcal{H}^{\mu}=2\frac{\sigma_{1}\sigma_{3}}{r_2r_4}(x^2+\xi^2)v.$$
	Using \eqref{coh2}, we get
	\begin{equation}\label{exp1}
		\mathrm{R}_{n}(x)=\mathrm{H}_{n}^{\mu}(x)+\varepsilon_{n-2}\mathrm{H}_{n-2}^{\mu}(x) , \ \  n\geq2.
	\end{equation}

On the other hand,
	$$\begin{cases}
	\mathrm{H}_{0}^{\mu}(x)=1, \ \ \ \mathrm{H}_{1}^{\mu}(x)=x,\\
		\mathrm{H}_{n+1}^{\mu}(x)=x\mathrm{H}_n^{\mu}(x)- \gamma_{n} \mathrm{H}_{n-1}^{\mu}(x), \ \ n \geq 1.
	\end{cases}$$
	with
	$$\gamma_{2n}=n,  \gamma_{2n-1}= n - 1/2  +\mu),  \\ \ n\geq 1,$$
Notice that the positive definite  condition yields $\mu > -n -\frac{1}{2}$, $n \geq 0$.
	
	We denote the TTRR relation for the MOPS $\{R_n\}_{n\geq0}$
	$$\mathrm{R}_{n+1}(x)=x\mathrm{R}_n(x)-\tilde{\gamma}_{n}\mathrm{R}_{n-1}(x),\ \ n\geq1, \ \ \mathrm{R}_{1}(x)=x, \ \ \mathrm{R}_{0}(x)=1.$$
	From Proposition \ref{prop},
	$$\tilde{\gamma}_{1}=\frac{1}{2}+\mu-\varepsilon_{0},$$
	$$\tilde{\gamma}_{n}=\frac{1}{2}\mu_n+\varepsilon_{n-2}-\varepsilon_{n-1}, \ \ n \geq 2,$$
	and
	$$\frac{1}{2}\mu_n\varepsilon_{n}=\varepsilon_{n-1}(\frac{1}{2}\mu_{n+2}+\varepsilon_{n}-\varepsilon_{n+1}), \ \ n \geq 1.$$
	In addition, by \eqref{eqdiff}, the sequence $\{\varepsilon_{n}\}_{n\geq0}$ satisfies the nonlinear quadratic difference equation
	$$\varepsilon_{n+1}=\frac{1}{2}(\mu_{n+1}+\mu_{n+2})+(\varepsilon_{1}-1)(1-\frac{1+2\mu}{2\varepsilon_{0}})-\frac{1}{4}\frac{\mu_{n}\mu_{n+1}}{\varepsilon_{n-1}}, \ \ n \geq 1.$$
	Therefore, for  fixed initial values $\varepsilon_{0}$ and $\varepsilon_{1}$ we can compute the sequences of parameters $\{\varepsilon_{n}\}_{n\geq0}$ and $\{\tilde{\delta}_{n}\}_{n\geq1}$.
\bigskip

As an illustrative example, if  $\mu=5$, $\lambda=0.1$, $\varepsilon_{0}=1.2$, $\varepsilon_{1}=1.3$ and $\xi=0,$ then
$\gamma_n=\frac{1}{2}\mu_n$, $s_{0}^{(\lambda,\mu)}=1$, $s_{1}^{(\lambda,\mu)}=5.6$, $s_{2}^{(\lambda,\mu)}=7.22$, $\eta_{-1}^\mu(\lambda)=0$, $\eta_{0}^\mu(\lambda)=1.39$, $\eta_{1}^\mu(\lambda)=1.98$.\\
Furthermore,
$$u=\mathcal{H}^5=0.02x^{10}e^{-x^2}.$$
As a consequence,
$$v=0.45x^8e^{-x^2}+\delta_0.$$

Next, let consider the function
$$f(x)=x(10-x).$$
Obviously $f \in \mathcal{W}_2^1(\mathbb{R}, u, v, \mu)$. In this case $f_{0}^{(\lambda, \mu)}=-5.758$, $f_{1}^{(\lambda, \mu)}=287.886$, $w_{0}^\mu(f)=-15.18$ and $w_{1}^\mu(f)=80.551$.

Now, applying the algorithm introduced in Section $6$, we find the Fourier $\mathcal{T}_{\mu}$–Sobolev coefficients of $f$.

\subsection{The generalized Gegenbauer case}
The sequence of generalized Gegenbauer polynomials $\{\mathrm{G}_n^{(\alpha,\mu-\frac{1}{2})}\}_{n\geq0}$  is an MOPS with respect to the positive definite linear form $\mathcal{G}^{(\alpha,\mu-\frac{1}{2})}$ defined by the weight function
	$$\omega(x)=\frac{\Gamma(\alpha+\mu+\frac{3}{2})}{\Gamma(\alpha+1)\Gamma(\mu+\frac{1}{2})}|x|^{2\mu}(1-x^2)^{\alpha}, \alpha >-1, $$
supported on the interval $ (-1,1).$

The corresponding form  is a positive definite $\mathcal{T}_{\mu}$-classical linear form and satisfies the  $\mathcal{T}_{\mu}$-Pearson equation
	$$\mathcal{T}_{\mu}((x^2-1)\mathcal{G}^{(\alpha,\mu-\frac{1}{2})})=2(\alpha+1)x\mathcal{G}^{(\alpha,\mu-\frac{1}{2})}.$$
	
The MOPS $\{\mathrm{G}_n^{(\alpha,\mu-\frac{1}{2})}\}_{n\geq0}$ satisfies the TTRR  \eqref{recurrence} (see \cite{Belmehdi2001})
	\begin{gather*}
		\hat{\beta}_n=0,\ \ \ \gamma_{n+1}=\frac{(n+1+\varrho_n)(n+1+2\alpha+\varrho_n)}{4(n+\alpha+\mu+\frac{1}{2})(n+\alpha+\mu+\frac{3}{2})},\ \ n\geq0,\\
		\varrho_n=2\mu\frac{1+(-1)^n}{2}, \ \ n\geq0.
	\end{gather*}	
If we apply $\mathcal{T}_{\mu}$ to $\mathrm{G}_n^{(\alpha,\mu-\frac{1}{2})}$, then we obtain \cite{bencheikh2007}
	\begin{equation}\label{geg}
		\mathcal{T}_{\mu}(\mathrm{G}_n^{(\alpha,\mu-\frac{1}{2})})=\mu_n\mathrm{G}_{n-1}^{(\alpha+1,\mu-\frac{1}{2})}.
	\end{equation}
	
	Let $v$ be a symmetric positive measure and denote by $\{R_n\}_{n\geq0}$ the corresponding MOPS. Assume that $(\mathcal{G}^{(\alpha,\mu-\frac{1}{2})} , v)$  is a $\mathcal{T}_{\mu}$-symmetric coherent pair. Thus, from \cite{sghaier2021}, we get
	$$(x^2-1)\mathcal{G}^{(\alpha, \mu-\frac{1}{2})}=2\frac{\sigma_{1}\sigma_{3}}{r_2r_4}(x^2-\xi^2)v, |\xi| \geq 1.$$
	Since
	$$(1- x^2)\mathcal{G}^{(\alpha, \mu-\frac{1}{2})}=\frac{\alpha+1}{\alpha+ \mu+\frac{3}{2}}\mathcal{G}^{(\alpha+1, \mu-\frac{1}{2})},$$
	then
	$$\mathcal{G}^{(\alpha+1, \mu-\frac{1}{2})}=2\frac{\alpha+ \mu+\frac{3}{2}}{\alpha+1}\frac{\sigma_{1}\sigma_{3}}{r_2r_4}(\xi^2- x^2)v, \xi \geq 1.$$
	
	Now, from  \eqref{coh2} and \eqref{geg}
	$$R_n(x)= \mathrm{G}_n^{(\alpha+1,\mu-\frac{1}{2})}(x)+\varepsilon_{n-2}\mathrm{G}_{n-2}^{(\alpha+1,\mu-\frac{1}{2})}(x),  n\geq2,$$
	where $\varepsilon_{n}$  is non-zero constant for all $n \geq 0$.

	The TTRR  for the MOPS $\{R_n\}_{n\geq0}$ reads as
	$$ \mathrm{R}_{n+1}(x)=x\mathrm{R}_{n}(x)-\tilde{\gamma}_{n} \mathrm{R}_{n-1}(x), \ \ n\geq1, \ \ \mathrm{R}_0(x)=1, \ \ 	\mathrm{R}_1(x)=x,$$
	with $\tilde{\gamma}_{n}>0$ for all $n\geq 1$. Then, by Corollary \ref{relation}
	$$\tilde{\gamma}_{n}=\frac{(n+\varrho_{n-1})(n+2\alpha+\varrho_{n-1})}{4(n+\alpha+\mu-\frac{1}{2})(n+\alpha+\mu+\frac{1}{2})}+\varepsilon_{n-2}-\varepsilon_{n-1}, \ \ n\geq 1.$$
	In addition, according to \eqref{eqdiff},  we can write the parameters $\{\varepsilon_n\}_{n\geq 0}$ as
	\begin{align*}
		\varepsilon_{n+1}& =\frac{1}{4(n+\alpha+\mu+\frac{3}{2})}\bigg(  \frac{(n+1+\varrho_n)(n+1+2\alpha+\varrho_n)}{n+\alpha+\mu+\frac{1}{2}} +\frac{(n+2+\varrho_{n+1})(n+2+2\alpha+\varrho_{n+1})}{n+\alpha+\mu+\frac{5}{2}}\bigg)\\
		& \qquad +\Lambda-\frac{(n+\varrho_{n-1})(n+2\alpha+\varrho_{n-1})(n+1+\varrho_n)(n+1+2\alpha+\varrho_n)}{16(n+\alpha+\mu+\frac{1}{2})^2(n+\alpha+\mu-\frac{1}{2})(n+\alpha+\mu+\frac{3}{2})\varepsilon_{n-1}},
	\end{align*}
	where
	$$\Lambda=\Big(\varepsilon_1-\frac{1+\alpha}{(\frac{3}{2}+\alpha+\mu)(\frac{5}{2}+\alpha+\mu)}\Big)\Big(1-\frac{\mu_1(1+2\alpha+2\mu)}{4\varepsilon_0(\frac{1}{2}+\alpha+\mu)(\frac{3}{2}+\alpha+\mu)}\Big).$$

	Therefore, for  fixed initial values $\varepsilon_{0}$ and $\varepsilon_{1}$, we can compute the sequences of parameters $\{\varepsilon_{n}\}_{n\geq0}$ and $\{\tilde{\gamma}_{n}\}_{n\geq1}$.	
	
	 We choose $\alpha=5$, $\lambda =0.5$, $\mu=1$,  $\varepsilon_{0}=0.1$,  $\varepsilon_{1}=0.15$ and $\xi=1$. In this case, $$\gamma_{n+1}=\frac{(n+1+\varrho_{n})(n+11+\varrho_{n})}{4(n+\frac{13}{2})(n+\frac{15}{2})}, \ \ \ \varrho_{n}=1+(-1)^n, \ \ n \geq 0,$$
 $\Lambda=0.0822$, $s_{0}^{(\lambda,\mu)}=1$, $s_{1}^{(\lambda,\mu)}=0.7$, $s_{2}^{(\lambda,\mu)}=0.219$, $\eta_{-1}^\mu(\lambda)=0$, $\eta_{0}^\mu(\lambda)=0.049$, $\eta_{1}^\mu(\lambda)=0.0258$.\\
 Furthermore
 $$u=\mathcal{G}^{(5,\frac{1}{2})}=17.6x^{2}(1-x^2)^{5},\ \ \ -1<x<1.$$
 Therefore
 $$v=-4.3325.10^{-6}(1-x^2)^{-1}\mathcal{G}^{(6,\frac{1}{2})}+\frac{1}{2}(\delta_1+\delta_{-1}),$$
 where
 $$\mathcal{G}^{(6,\frac{1}{2})}=22x^{2}(1-x^2)^{6},\ \ \ -1<x<1.$$
 So
 $$v=-0.953.10^{-4}x^{2}(1-x^2)^5
  + \frac{1}{2}(\delta_1+\delta_{-1}),\ \ \ -1<x<1.$$
 We will use the function
 $$f(x) = xe^{-(x-0.2)^2}.$$
 It is clear that $f \in \mathcal{W}_2^1(\mathbb{R}, u, v, \mu)$.
  Then $f_{0}^{(\lambda, \mu)}=0.0583$, $f_{1}^{(\lambda, \mu)}=1.338$, $w_{0}^\mu(f)=0.081$ and $w_{1}^\mu(f)=1.768$.

  These initial conditions and the algorithm introduced in the previous section yield the Fourier Dunkl–Sobolev coefficients of $f$.
  \bigskip
  \bigskip


 \section{Acknowledgements}
 The work of Francisco Marcell\'an has
been supported by the research project [PID2021- 122154NB-I00], \emph{Ortogonalidad y Aproximación con Aplicaciones en Machine Learning y Teoría de la Probabilidad} funded  by MICIU/AEI/10.13039/501100011033 and by "ERDF A Way of making Europe”and the Madrid Government (Comunidad de Madrid-Spain) under the Multiannual Agreement with UC3M in the line of Excellence of University Professors, grant EPUC3M23 in the
context of the V PRICIT (Regional Program of Research and Technological Innovation).

	

\begin{thebibliography}{99}

\bibitem{Belmehdi2001} S. Belmehdi,  \emph{Generalized Gegenbauer orthogonal polynomials.} Proceedings of the Fifth International Symposium on Orthogonal Polynomials, Special Functions and their Applications (Patras, 1999)
J. Comput. Appl. Math. \textbf{133} (2001), no. 1-2, 195-205.

\bibitem{bencheikh2007} Y. Ben Cheikh, M. Gaied,\emph{ Characterization of the Dunkl-classical symmetric orthogonal polynomials.} Appl. Math. Comput. \textbf{187 }(2007), 105-114.


\bibitem{Cleonice-10} C. F. Bracciali, L.  Castaño-García, J. J. Moreno-Balc\'azar, \emph{Some asymptotics for Sobolev orthogonal polynomials involving Gegenbauer weights.} J. Comput. Appl. Math. \textbf{235 }(2010), no. 4, 904-915.

\bibitem{Laura-03} L. Castaño-García, J. J.  Moreno-Balc\'azar,  \emph{A Mehler-Heine-type formula for Hermite-Sobolev orthogonal polynomials.} J. Comput. Appl. Math. \textbf{150} (2003), no. 1, 25-35.

\bibitem{ted} T. S. Chihara, \emph{Generalized Hermite polynomials.} Doctoral Dissertation, Purdue, 1955.

\bibitem{chihara1978} T. S. Chihara,\emph{ An Introduction to Orthogonal Polynomials.} Gordon and Breach, New York, 1978.	

\bibitem{Oscar-SIGMA2018}  O. Ciaurri, J.  M\'inguez,  \emph{Fourier series of Gegenbauer-Sobolev polynomials.} SIGMA Symmetry Integrability Geom. Methods Appl. \textbf{14} (2018), Paper no. 024, 11 pp.

\bibitem{Marcel-07} M. G. De Bruin, W. G.  M. Groenevelt, F. Marcell\'an, H. G.  Meijer, J. J.  Moreno-Balc\'azar,  \emph{Asymptotics and zeros of symmetrically coherent pairs of Hermite type.}  In \emph{Difference Equations, Special Functions and Orthogonal Polynomials}, 378-393, World Sci. Publ., Hackensack, NJ, 2007.

\bibitem{Delgado2005} A. M. Delgado, F. Marcell\'an, \emph{On an extension of symmetric coherent pairs of orthogonal polynomials.}  J. Comput. Appl. Math.  \textbf{178} (2005),  no. 1-2, 155--168.

\bibitem{herbert} H. Due\~nas Ruiz, F. Marcell\'an, A. Molano, \emph{On symmetric (1,1)-coherent pairs and Sobolev orthogonal polynomials: an algorithm to compute Fourier coefficients}. Rev. Colombiana Mat. \textbf{53} (2019), no. 2, 139-164.

\bibitem{dunkl1991} C. F. Dunkl, \emph{ Integral kernels reflection group invariance.} Canad. J. Math. \textbf{43} (1991), 1213-1227.

\bibitem{Iserles1991} A. Iserles, P. E. Koch, S. P. Nørsett, J. M. Sanz-Serna, \emph{On polynomials orthogonal with respect to certain Sobolev inner products.} J. Approx. Theory \textbf{65}  (1991), no. 2,  151–175.

\bibitem{marcellan1994} F. Marcell\'{a}n, A. Branquinho, J. Petronilho, \emph{Classical orthogonal polynomials: A functional approach.} Acta Appl. Math. \textbf{34} (1994), 283-303.

\bibitem{maite} F. Marcell\'an, T. E.  P\'erez, M. A. Pi\~nar,  \emph{Gegenbauer-Sobolev orthogonal polynomials.} Math. Appl.,\textbf{ 296}. Kluwer Academic Publishers Group, Dordrecht, 1994, 71–82.

\bibitem{zecarlos}  F. Marcell\'an, J. Petronilho, \emph{Orthogonal polynomials and coherent pairs: the classical case.} Indag. Math. (N. S.) \textbf{6} (1995), no. 3, 287-307.

\bibitem {Paco-2006} F. Marcell\'an, J. J. Moreno-Balc\'azar,  \emph{Asymptotics and zeros of Sobolev orthogonal polynomials on unbounded supports.} Acta Appl. Math.  \textbf{94} (2006), no. 2, 163-192.

\bibitem{Paco2015} F. Marcell\'an, Yuan  Xu,  \emph{On Sobolev orthogonal polynomials.} Expo. Math. \textbf{33} (2015), no. 3, 308-352.

\bibitem{pascal} P. Maroni, \emph{Sur la suite de polyn\^{o}mes orthogonaux associ\'ee \`a la forme $u=\delta_{c} + \lambda (x-c)^{-1} L.$} Period. Math. Hungar. \textbf{21} (1990), no. 3, 223-248.

\bibitem{andrei1997} A.  Mart\'inez-Finkelshtein, J. J. Moreno-Balc\'azar, H. Pijeira-Cabrera,\emph{ Strong asymptotics for Gegenbauer-Sobolev orthogonal polynomials.} J. Comput. Appl. Math. \textbf{81} (1997), no. 2, 211-216.

\bibitem{Meijer1997} H. G. Meijer,  \emph{Determination of all coherent pairs.} J. Approx. Theory \textbf{89} (1997), no. 3, 321--343.

\bibitem{Ros1994} M. Rosenblum, \emph{Generalized Hermite polynomials and the Bose-like oscillator calculus.} Oper. Theory Adv. Appl. \textbf{73} (1994), 369-396.	

\bibitem{sghaier2012}	 M. Sghaier  \emph{A note on the Dunkl-classical orthogonal polynomials.} Integral Transforms Spec. Funct. \textbf{23} (2012), 753–760.


\bibitem{sghaier2021}  M. Sghaier, S. Hamdi, \emph{Some results on Dunkl-coherent pairs.} Integral Transforms Spec. Funct. \textbf{32} (2021), no. 11, 863-882.

\bibitem{alexei} A. Zhedanov,\emph{ Rational spectral transformations and orthogonal polynomials.} J. Comput. Appl. Math. \textbf{85} (1997), 67-86.
	
	
		
		
		
		
	
		
		
	 	
	\end{thebibliography}
\end{document}